\documentclass{amsart}
\usepackage[shortlabels]{enumitem}

\numberwithin{equation}{section}

\usepackage{amssymb,amsmath,amsthm,verbatim,IEEEtrantools,hyperref}
\theoremstyle{definition}

\newtheorem*{mythm*}{Theorem}
\newtheorem*{mydef*}{Definition}
\newtheorem*{myex*}{Example}
\newtheorem*{mycor*}{Corollary}
\newtheorem*{mypro*}{Proposition}
\newtheorem*{mylem*}{Lemma}
\newtheorem*{mycon*}{Conjecture}
\newtheorem*{myprb*}{Problem}
\newtheorem*{myrmk*}{Remark}

\theoremstyle{definition}
\newtheorem{theorem}{Theorem}[section]
\newtheorem{thm}[theorem]{Theorem}
\newtheorem{lemma}[theorem]{Lemma}

\newtheorem{prop}[theorem]{Proposition}
\newtheorem{example}[theorem]{Example}
\newtheorem{remark}[theorem]{Remark}

\def\RR{\mathbb{R}}
\def\CC{\mathbb{C}}

\def\QQ{\mathbb{Q}}
\def\ZZ{\mathbb{Z}}
\def\NN{\mathbb{N}}

\def\supp{\mathrm{supp}}

\newcommand{\rbr}[1]{\left( {#1} \right)}

\newcommand{\cbr}[1]{ \left\{ {#1} \right\} }
\newcommand{\abr}[1]{\left\langle {#1} \right\rangle}
\newcommand{\abs}[1]{\left| {#1} \right|}

\begin{document}

\author{Robert Fraser \and Kyle Hambrook}
\title{Explicit Salem Sets in $\RR^n$}
\begin{abstract}
We construct the first explicit (i.e., non-random) examples of 
Salem sets in $\mathbb{R}^n$ of arbitrary prescribed Hausdorff dimension.  
This completely resolves a problem proposed by Kahane more than 60 years ago. 
The construction is based on a form of Diophantine approximation in number fields. 
\end{abstract}

\maketitle

\section{Main Results}

%\textbf{Notable Notation.} 
%For $x \in \RR$, $e(x)=e^{-2\pi i x}$. 
%
%For $x \in \RR^n$, $|x|$ is the max-norm of $x$, and $|x|_2$ is the 2-norm of $x$; 
%i.e., $|x| = \max_{1 \leq i \leq n} |x_i|$ and $|x|_2 = \rbr{ \sum_{i=1}^{n} |x_i|^2 }^{1/2}$. 
%
%For $x \in \RR^n$, $|x| = \max_{1 \leq i \leq n} |x_i|$ is the max-norm of $x$, 
%and $|x|_2 = \rbr{ \sum_{i=1}^{n} |x_i|^2 }^{1/2}$ is the 2-norm of $x$. 
%

For $x \in \RR^n$, $|x| = \max_{1 \leq i \leq n} |x_i|$, i.e., the max-norm of $x$. 
%For $x \in \RR$, $e(x)=e^{-2\pi i x}$. 
%
Let $K$ be a number field (i.e., a finite extension field of $\QQ$) of degree $n$. 
%Fix a number field $K$ of degree $n$. 
Let $B=\cbr{\omega_1,\ldots,\omega_{n}}$ be an integral basis for $K$. 
We identify $\mathbb{Q}^n$ with $K$ by identifying 
$q=(q_1,\ldots,q_n) \in \mathbb{Q}^n$ with 
$q=\sum_{i=1}^{n} q_i \omega_i \in K$. 
Since $B$ is an integral basis, 
this also identifies $\mathbb{Z}^n$ with $\mathcal{O(K)}$, 
the ring of integers for $K$. 
%By extension, we identify $\RR^n$ with the algebra $\RR \omega_1 + \cdots + \RR \omega_n$. 
%This makes $\RR^n$ into an algebra, $\QQ^n$ into a field, and $\ZZ^n$ into a ring. 
%
Let $\tau > 1$. 
%Define $E(K,B,\tau)$ to be the set of all $x \in \RR^n$ such that 
%$$
%|x-q^{-1}r| \leq |q|^{-(\tau+1)}
%$$
%for infinitely many $(q,r) \in \ZZ^n \times \ZZ^n$. 
Define 
$$
E(K,B,\tau) = \cbr{x \in \RR^n : |x-r/q| \leq |q|^{-(\tau+1)} \text{ for infinitely many } (q,r) \in \ZZ^n \times \ZZ^n}.
$$

Our main result is the following theorem. 

\begin{thm}\label{main thm}
$E(K,B,\tau)$ is a Salem set of dimension $2n/(1+\tau)$. 
\end{thm}

Since $\tau > 1$ is arbitrary, 
Theorem \ref{main thm} yields Salem sets of every dimension 
$s \in (0,n)$. 
(For the endpoints, note that $\emptyset$ and $\RR^n$ 
are trivial examples of Salem sets of dimension $0$ and $n$, 
respectively.)  

%In the course of proving this theorem, we will also prove the following theorem 
%The following theorem is the main part of the proof of the theorem above. 

%The main part of the proof of Theorem \ref{main thm} is the following theorem. 

%As will be explained in Section \ref{background and motivation}, 
%Theorem \label{main thm} is saying that the Hausdorff and Fourier dimensions of $E(K,B,\tau)$ 
%are both equal to $2n/(1+\tau)$. 
%A simple covering argument shows that the 
%Hausdorff dimension of $E(K,B,\tau)$ is $\leq 2n/(1+\tau)$ 
%(Proposition \ref{hausdorff upper bound}). 
%The following existence theorem 

As will be evident from the definitions in Section \ref{background and motivation}, 
Theorem \ref{main thm} follows immediately from 
a simple upper bound on the Hausdorff dimension of $E(K,B,\tau)$ 
(namely, Proposition \ref{hausdorff upper bound}) 
and the following existence theorem, 
which is our main technical result. 

\begin{thm}\label{main-thm-2}
Let $r_1$ be the number of real embeddings of $K$ 
into $\CC$ and let 
and $r_2$ be the number of conjugate pairs of 
complex embeddings of $K$ into $\CC$.  
There exists a Borel probability measure $\mu$ 
with compact support contained in $E(K,B,\tau)$ such that 
\begin{align*}
\widehat{\mu}(\xi) 
= o
\rbr{
|\xi|^{-n/(1+\tau)} 
\exp \rbr{ \frac{ n \log |\xi| }{ \log \log |\xi| } } 
\log^{r_1 + r_2 + 1}|\xi| 
}
\quad \text{as $|\xi| \to \infty$}. 
\end{align*}
%for $\xi \in \RR^n$ with $|\xi| \geq 3$. 
%as $|\xi| \to \infty$. 
\end{thm}

%\textbf{Outline.} TBA

\noindent \textbf{Notation.} 
The expression $X \lesssim Y$ means $X \leq CY$ for some positive constant $C$ whose precise value is immaterial in the context. 
The expression $X \lesssim_{\alpha} Y$ has the same meaning, except the constant $C$ is permitted to depend also on a parameter $\alpha$. 
The expression $X \gtrsim Y$ means $Y \lesssim X$. 
The expression $X \approx Y$ means both $X \lesssim Y$ and $Y \lesssim X$.

\section{Background and Motivation}
\label{background and motivation}

For $x \in \RR$, $e(x)=e^{-2\pi i x}$. 
If $\mu$ is a finite Borel measure on $\mathbb{R}^n$, 
then the Fourier transform of $\mu$ is defined by
$$
\widehat{\mu}(\xi) 
= \int_{\mathbb{R}^n} e(x \cdot \xi) d\mu(x) 
\quad \text{ for all } \xi \in \mathbb{R}^n. 
$$

%Frostman's lemma \cite{Frostman} implies 
%that 
The Hausdorff dimension  
$\dim_H(E)$ of a Borel set $E \subseteq \mathbb{R}^n$ 
is equal to the supremum of the values of $s \in [0,n]$ such that the integral 
\[\int_{\RR^n} |\hat \mu(\xi)|^2 \, |\xi|^{s-n}d \xi \]
is convergent for some probability measure $\mu$ supported on $E$. 
This characterization of Hausdorff dimension is well-known; 
see for example 
\cite{falconer-book-1}, \cite{federer}, \cite{mattila-book-1}, \cite{mattila-book-2}, \cite{wolff-book}. 
It can be viewed as the statement that, for any $\epsilon > 0$, 
the Fourier transform $|\hat \mu(\xi)|$ of $\mu$ decays like 
$|\xi|^{-s/2 + \epsilon}$ in $L^2$-average. 
%an $L^2$-average sense. 

In contrast, the Fourier dimension $\dim_F(E)$ of a set $E \subseteq \mathbb{R}^n$ concerns the fastest 
pointwise rate of decay of the Fourier transform. 
The Fourier dimension of a set $E \subseteq \mathbb{R}^n$ is defined to be the 
supremum of the values of $s \in [0,n]$ such that
\[
|\hat \mu(\xi)|^2 |\xi|^{-s} \to 0 \quad \text{as $|\xi| \to \infty$}
\]
for some probability measure $\mu$ supported on $E$.

As general references for Hausdorff and Fourier dimension, 
see \cite{falconer-book-1}, \cite{federer}, \cite{mattila-book-1}, \cite{mattila-book-2}, \cite{wolff-book}. 
Recent papers by Ekstr\"{o}m, Persson, and Schmeling \cite{EPS} and Fraser, Orponen, and Sahlsten \cite{FOS} 
have revealed some interesting subtleties about Fourier dimension. 

Immediately from the definitions, we see that for every Borel set $E \subseteq \mathbb{R}^n$, 
$$
\dim_F(E) \leq \dim_H(E).
$$
%From this inequality and the definition of Fourier dimension, it is evident that 
%Theorem \ref{main thm} follows from Theorem \ref{main-thm-2} 
%and the Hausdorff dimension upper bound furnished by Proposition \ref{hausdorff upper bound}. 
%Moreover, it is clear that Theorem \ref{main thm} follows from Theorem \ref{main-thm-2} 
%and the Hausdorff dimension upper bound furnished by Proposition \ref{hausdorff upper bound}. 

Every $k$-dimensional plane in $\mathbb{R}^n$ with $k < n$ has Fourier dimension $0$ and Hausdorff dimension $k$. 
More generally, every subset of every $(n-1)$-dimensional plane in $\RR^n$ has Fourier dimension $0$, while the Hausdorff dimension may take any value in $[0,n-1]$. 
%If $A$ is a $k$-dimensional subspace in $\mathbb{R}^d$ with $k < d$, then $\dim_F A = 0$ and $\dim_H A = k$. 
%If $A$ is a $k$-dimensional subspace in $\mathbb{R}^d$ with $k < d$, then $\dim_F A = 0$ and $\dim_H = k$. 
The middle-$1/3$ Cantor set in $\mathbb{R}$ has Fourier dimension $0$ and Hausdorff dimension $\log 2 / \log 3$. 
%K{\"o}rner \cite{Korner} has shown that for every $0 \leq s \leq t \leq 1$ there is a compact set $A \subseteq \mathbb{R}$ with Fourier dimension $s$ and Hausdorff dimension $t$. 
More generally, a middle-$\delta$ Cantor set in $\mathbb{R}$ may have positive Fourier dimension;
however, its Fourier dimension will always be strictly smaller than its Hausdorff dimension. 

Sets $E \subseteq \mathbb{R}^n$ with $$\dim_F(E) = \dim_H(E)$$ are called Salem sets. 

%A set $A \subseteq \mathbb{R}^d$ with $\dim_F A = \dim_H A$ is called a Salem set. 
%Every ball in $\mathbb{R}^d$ is a Salem set of dimension $d$. 
Every set in $\RR^n$ that contains a ball is a Salem set of dimension $n$. 
%Every countable set in $\mathbb{R}^d$ is a Salem set of dimension zero. 
Every set in $\RR^n$ of Hausdorff dimension $0$ is a Salem set of dimension $0$. 
Less trivially, every sphere in $\mathbb{R}^n$ 
(or, more generally, every $(n-1)$-dimensional manifold in $\RR^n$ 
with non-vanishing Gaussian curvature) 
is a Salem set of dimension $n-1$. 

%Salem sets in $\mathbb{R}^d$ of dimension 
%%%$s \notin \cbr{0,d-1,d}$   
%$s \neq 0,d-1,d$   
%are harder to find.
%
Salem sets in $\mathbb{R}^n$ of dimension 
%$s \notin \cbr{0,d-1,d}$   
$s \neq 0,n-1,n$   
%are more exotic. 
%are harder to find.
are more complicated. 

It is known that given any $s \leq n$, 
there exist Salem sets of dimension $s$ contained in $\mathbb{R}^n$. 
%There are many random constructions of Salem sets. 
Using Cantor sets with randomly chosen contraction ratios, 
Salem \cite{Salem} was the first to show that for every $s \in [0,1]$ 
there is a Salem set in $\mathbb{R}$ of dimension $s$. 
Kahane showed that images of compact subsets of $\mathbb{R}^d$ 
under certain stochastic processes 
(namely, Brownian motion, fractional Brownian motion, and Gaussian Fourier series) 
are almost surely Salem sets 
(see {\cite{Kahane-1966-Brownian}, \cite{Kahane-1966-Fourier}, 
\cite[Ch.17,18]{kahane-book}). 
%As a consequence of these results, 
Through these results, 
Kahane established that for every $s \in [0,n]$ 
there is a Salem set in $\mathbb{R}^n$ of dimension $s$.  
%Ekstr\"{o}m \cite{ek-2} has shown that every Borel set in $\mathbb{R}$ is diffeomorphic to a Salem set in $\mathbb{R}$ under a random diffeomorphism.  
Ekstr\"{o}m \cite{ekstrom} has showed that the image of 
any Borel set in $\mathbb{R}$ under a random diffeomorphism 
is almost surely a Salem set. 
Other random constructions of Salem sets have been given by 
Bluhm \cite{Bluhm-1}, {\L}aba and Pramanik \cite{LP}, 
Shmerkin and Suomala \cite{shmerkin-suomala}, 
and Chen and Seeger \cite{chen-seeger}.

%These random constructions give collections of sets where each individual set is ``almost surely'' or ``with positive probability" a Salem set. But they don't provide any explicit examples of Salem sets. 

%However, most of these Salem set constructions are probabilistic in nature, 
%revolving around random cantor sets 
%\cite{Salem51}, Brownian motion \cite{Kahane66}, \cite{Kahane66B}, 
%or some other random procedure. 

Kahane \cite{kahane-book} suggested that it would be interesting to find 
\textit{explicit} (by which Kahane meant \textit{non-random}) constructions of 
Salem sets in $\RR^n$ of every dimension $s \in [0,n]$. 

Explicit Salem sets of dimensions $0$, $n-1$, or $n$ are easy to find. 
Indeed, see the examples of Salem sets of dimensions $0$, $n-1$, or $n$ we listed above. 

All known explicit examples of Salem sets of dimension other than 
$0, n-1$ or $n$ in $\mathbb{R}^n$ are based on a construction by 
Kaufman \cite{Kaufman81}. 
Kaufman considered sets of numbers that are well-approximated by real numbers.
%Given any $\tau>1$, Kaufman defines the set $E(\tau)$ to be the set
%\[
%\{x \in [0,1] : |qx - r| \leq |q|^{-\tau} \quad 
%\text{for infinitely many pairs of integers $(q,r)$}\}
%\]
For $\tau>1$, Kaufman studied the set 
\[
E(\tau) = \cbr{x \in \mathbb{R} : |xq - r| \leq |q|^{-\tau} 
\text{ for infinitely many } (q,r) \in \ZZ \times \ZZ}.
\] 
Much earlier, Jarn\'{\i }k \cite{Jarnik29} and 
Besicovitch \cite{Besicovitch34} 
showed that for $\tau > 1$, 
the set $E(\tau)$ has Hausdorff dimension 
equal to ${2}/{(1 + \tau)}$. 
This is a key result in metric Diophantine approximation. 
Kaufman \cite{Kaufman81} established 
pointwise Fourier decay bounds 
for a natural measure supported on the set $E(\tau)$, 
thereby showing that the Fourier dimension of 
$E(\tau)$ is also equal to 
${2}/{(1 + \tau)}$ for $\tau > 1$. 
This provides explicit Salem sets in $
\RR$ of arbitrary dimension $s \in (0,1)$.  
(Note that Dirichlet's approximation theorem 
gives $E(\tau) = \mathbb{R}$ 
when $\tau \leq 1$.)  

%Bluhm \cite{Bluhm-2} 
%gave a detailed account of what is essentially 
%Kaufman's proof and also 
%pointed out that 
%(as a consequence of a theorem of Gatesoupe \cite{Gatesoupe}) 
Bluhm combined Kaufman's argument 
with a theorem of Gatesoupe \cite{Gatesoupe} 
to show that the rotationally symmetric set 
$$
\cbr{x \in \mathbb{R}^n : |x|_2 \in E(\tau)} 
$$ 
(where $|x|_2 = ( \sum_{i=1}^{n} |x_i|^2 )^{1/2}$ 
is the 2-norm) 
is 
%an explicit Salem set
a Salem set 
in $\mathbb{R}^n$ of dimension $n-1 + 2/(1+\tau)$ 
whenever $\tau > 1$. 
This gives explicit Salem sets in $\RR^n$ of every dimension $s \in (n-1,n)$, 
but leaves open the range $s \in (0,n-1)$. 
%However, explicit Salem sets in $\mathbb{R}^d$ of dimension 
%$0 < s < d-1$ were unknown until now.

In metric Diophantine approximation, 
the natural multi-dimensional generalization of $E(\tau)$ is 
%$$
%E(m,n,\tau) = \cbr{x \in \mathbb{R}^{mn} : \max_{1 \leq i \leq m}|\sum_{j=1}^{n} q_j x_{ij} %- r_i| \leq |q|^{-\tau} \text{ for infinitely many } (q,r) \in \ZZ^n \times \ZZ^m},
%$$
$$
E(m,n,\tau) = \cbr{x \in \mathbb{R}^{mn} : |xq - r| \leq |q|^{-\tau} 
\text{ for infinitely many } (q,r) \in \ZZ^n \times \ZZ^m},
$$
where we identify $\mathbb{R}^{mn}$ with the set of $m \times n$ matrices 
with real entries, 
so that $xq$ is computed as the product of an $m \times n$ 
and an $n \times 1$ column vector. 
By Minkowski's theorem on linear forms, $E(m,n,\tau) = \mathbb{R}^{mn}$ 
when $\tau \leq n/m$.  
%Jarn{\'\i}k \cite{Jarnik} and Eggleston \cite{Eggleston} showed 
%the Hausdorff dimension of $E(m,1,\tau)$ is $(m+1)/(1+\tau)$ if $\tau > 1/m$. 
Bovey and Dodson \cite{bovey-dodson} proved that the Hausdorff dimension of 
$E(m,n,\tau)$ is $m(n-1) + (m+n)/(1+\tau)$ when $\tau > n/m$. 
The $n=1$ case was established earlier by Jarn{\'\i}k \cite{jarnik-1} 
and Eggleston \cite{Eggleston}. 
The mass transference principle and slicing technique of Beresnevich 
and Velani \cite{BV-annals}, \cite{BV-slicing} 
may also be used to compute the Hausdorff dimension of $E(m,n,\tau)$.  
%The second author 
Hambrook 
\cite{hambrook-tams} proved that the Fourier dimension of 
$E(m,n,\tau)$ is at least $2n/(1+\tau)$ if $\tau > n/m$. 
However, there is a gap between the Hausdorff dimension and 
this lower bound on the Fourier dimension, 
and so it not known  
%However, it is unclear 
whether $E(m,n,\tau)$ is a Salem set when 
$\tau > n/m$ and $mn > 1$. 

The first explicit examples of Salem sets of arbitrary dimension in $\mathbb{R}^2$ 
are due to 
%the second author  
Hambrook \cite{Hambrook17}. 
The construction 
%relies on an 
uses an analogue of the set $E(\tau)$ described above. 
%An analogue of $E(\tau)$ is defined where, 
Instead of considering real numbers $x$ such that $x$ is close to many rational numbers 
$\frac{r}{q}$, one considers real vectors $(x_1, x_2)$ 
such that $x_1 + x_2 i$ is close to many ratios of Gaussian integers, 
i.e., close to many complex numbers of the form 
$\frac{r_1 + r_2 i}{q_1 + q_2 i}$, where $r_1, r_2, q_1, q_2$ are integers.  
Precisely, \cite{Hambrook17} 
%considers the set 
shows that the set 
$$
E(\CC, \tau) = \cbr{x \in \mathbb{R}^2 : |qx - r| \leq |q|^{-\tau} \text{ for infinitely many } (q,r) \in \ZZ^2 \times \ZZ^2}
$$
is Salem with dimension $4/(1+\tau)$ when $\tau > 1$. 
Here $\RR^2$ and $\CC$ are identified in the usual way, so $qx$ is viewed as a product of complex numbers. 
%Of course, this identifies $\ZZ^2$ with the Gaussian integers $\ZZ[i] =\ZZ + i \ZZ$. 
Of course, this identifies $\ZZ^2$ with the Gaussian integers $\ZZ[i]$ and $\QQ^2$ with the number field $\QQ(i)$. 

\section{Innovations}\label{Innovations}

The construction of explicit Salem sets in the present paper is inspired by the construction of Hambrook \cite{Hambrook17}, which is in turn inspired by the construction of Kaufman \cite{Kaufman81}. 
%In what follows, we describe our key innovations. 
We describe our key innovation in the following sequence of remarks. 
%the major points of interest along the path we took from the construction of \cite{Hambrook17} to the present construction. 

%\begin{enumerate}[(1)]

%\item 
\begin{remark}
The identification of 
$\RR^2$ and $\CC$ in \cite{Hambrook17} 
suggests identifying $\RR^n$ 
with some other algebraic structure 
and mimicking the argument. 
An reasonable idea is to identify $\RR^4$ 
with the set of quaternions. 
But, as explained in \cite{Hambrook17}, 
this does not seem to work. 
Our breakthrough idea was to 
%Instead, we 
shift focus to the subset $\QQ^n \subseteq \RR^n$. 
We identify $\QQ^n$ with a number field $K$ via 
an integral basis $B = \cbr{\omega_1,\ldots,\omega_n}$. 
Of course, this induces a identification of $\RR^n$ with 
the algebra $\RR\omega_1 + \cdots \RR\omega_n$, 
but we never use this. 
\end{remark}

%\item 
\begin{remark}
The inequality $|x-r/q| \leq |q|^{-(\tau+1)}$ that defines $E(K,B,\tau)$ 
is different from the inequality $|qx-r| \leq |q|^{-\tau}$ 
that defines $E(\CC,\tau)$. % and $E(\tau)$. 
(In the case of $E(\tau)$, the inequalities are actually equivalent 
because $|xy|=|x||y|$ for all $x,y \in \RR$.) 
Due to this difference in form, 
the estimation of a complex exponential sum naturally appears 
in our proof, whereas a complex exponential integral 
naturally appears in the proofs of Hambrook \cite{Hambrook17} 
and (implicitly) of Kaufman \cite{Kaufman81}. 
While it is possible to modify our proof to handle the version of the set 
$E(K,B,\tau)$ defined via the inequality $|qx-r| \leq |q|^{-\tau}$, 
we found the proof is easier with $E(K,B,\tau)$ as currently defined. 
Moreover, the proofs of Kaufman \cite{Kaufman81} for $E(\tau)$ 
and of Hambrook \cite{Hambrook17} for $E(\CC,\tau)$ (defined using either inequality)
can be modified to go through the complex exponential sum rather 
than the complex exponential integral. 
The method of proof via the complex exponential sum 
actually comes from \cite{hambrook-fraser-p-adic}, 
where a construction of explicit Salem sets in the $p$-adic numbers is given. 
\end{remark}

%\item 
\begin{remark} 
The proof in \cite{Hambrook17} 
requires the evaluation of a certain complex exponential integral 
(see Lemma 5). 
This evaluation involves only a simple calculation 
with dot products and real and imaginary parts. 
However, as we came to realize, this evaluation uses implicitly 
%
%In order to evaluate a certain complex exponential integral, 
%\cite[Lemma 5]{Hambrook17} used implicitly 
%
the property 
that the transpose of the usual matrix representation of an element of 
%$\CC$ 
$\QQ(i)$ 
%(or $\CC$)
coincides with the matrix representation of the complex conjugate 
of that element. 
Similarly, our proof requires the estimation 
of a certain complex exponential sum.  
%(see (2) above for a brief discussion of why 
%it is a sum rather than an integral). 
And, likewise, this sum turns out to depend on transposes of  
matrix representations of elements of $K$. 
Unfortunately, the analogous property about  
transposes of matrix representations  
does not hold in general number fields. 
To carry out the estimation, 
we require a 
%substantially 
more sophisticated 
property of transposes of matrix representations 
of elements of number fields. 
The details make up Section 
\ref{Algebraic Number Theory : An Exponential Sum section}. 
\end{remark}

%\item 
\begin{remark}
Kaufman's proof \cite{Kaufman81} depends on a simple prime divisor bound 
based on the uniqueness of prime factorization in $\ZZ$. 
The proof in \cite{Hambrook17} relies on the standard divisor bound in $\ZZ[i]$, 
which holds because $\ZZ[i]$ is a unique factorization domain with a finite unit group. 
However, the analogous divisor bound does not necessarily hold in the ring of integers $\ZZ_K$ 
of a general number field $K$. 
In general, 
%the ring of integers of a number field 
$\ZZ_K$
may not be a unique factorization domain,  
and its unit group may be infinite. 
In the present paper, 
we use the unique factorization of ideals 
and the geometric structure of the unit group of $\ZZ_K$ 
to prove a substitute divisor bound that turns out to be sufficient for our purpose. 
This substitute divisor bound and its proof is essentially due to Elkies \cite{ElkiesMO11}. 
The details make up Section \ref{Algebraic Number Theory: A Divisor Bound section}.  
\end{remark}

%\item 
\begin{remark}
In Section \ref{The Sets $Q(M)$, $Q'(M)$, and $Q''(M)$ section}, 
our proof requires the successive construction  of three sets: $Q(M) \supseteq Q'(M) \supseteq Q''(M)$. 
In short, $Q(M) = \cbr{q \in \ZZ^n : M/2 < |q|  \leq M}$,   
$Q'(M)$ is a subset of $Q(M)$ formed by removing those $q$ 
which are divisors (in $\ZZ_K$) 
of certain ``small'' elements of $\ZZ^n$, 
and  
$Q''(M)$ is a large subset of $Q'(M)$ 
whose elements $q$ have ideal norms $N(\abr{q})$ 
(defined in Section \ref{Algebraic Number Theory: Basics section}) 
that are all roughly the same size. 
In contrast, the proofs of Hambrook \cite{Hambrook17} 
and Kaufman \cite{Kaufman81} work (essentially)
with just the set $Q(M)$. 
The reason for the additional 
complication in our proof is ultimately is that, for elements 
$q$ in $\ZZ$ or $\ZZ^2 \approx \ZZ[i]$,  
the ideal norm $N(\abr{q})$ equals 
(respectively) $|q|$ or $|q|_2^2 = q_1^2+q_2^2$, 
while for elements $q \in \ZZ^n \approx \ZZ_K$ the ideal norm $N(\abr{q})$ 
is generally not comparable to $\|q\|^n$, 
where $\| {} \cdot {} \|$ is any norm on $\RR^n$. 
To be more specific, while it is true that $N(\abr{q}) \lesssim \|q\|^n$ for all $q \in \ZZ^n$, 
it is not true, for an arbitrary number field $K$ and basis $B$, that $N(\abr{q}) \gtrsim \|q\|^n$ for all $q \in \ZZ^n.$  
\end{remark}

%\end{enumerate}

\section{Hausdorff Dimension Upper Bound}\label{Hausdorff Dimension Upper Bound section}

%The required upper bound on the Hausdorff dimension of $E(K,B,\tau)$ is established by a standard covering argument. 

\begin{prop}\label{hausdorff upper bound}
$\dim_H(E(K,B,\tau)) \leq 2n/(1+\tau)$.  
\end{prop}
\begin{proof}
Note that $E(K,B,\tau)$ is invariant under translation by elements of $\ZZ^n$. 
Thus it suffices show that $E(K,B,\tau) \cap [-1/2,1/2]^n$ has Hausdorff dimension at most $2n/(1+\tau)$. 
Let $\overline{B}(x,r)$ denote the closed ball in $\RR^n$ with center $x$ and radius $r$. 
Note 
\begin{align*}
E(K,B,\tau) 
%&
%= \cbr{x \in \RR^n : \abs{x - \frac{r}{q}}  \leq |q|^{-(1+\tau)} 
%\text{ for $\infty$-many } (q,r) \in \ZZ^n \times \ZZ^n  } 
%\\
% &
 = \bigcap_{N=1}^{\infty} \bigcup_{|q| > N} \bigcup_{r \in \ZZ^n} 
 \overline{B}( r/q,|q|^{-(1+\tau)} ). 
\end{align*}
So, for every $N \geq 1$,  
\begin{align*}
E(K,B,\tau) \cap [-1/2,1/2]^n 
\subseteq \bigcup_{\substack{ q \in \ZZ^n \\ |q| > N}} 
\bigcup_{\substack{r \in \ZZ^n \\ |r| \leq |q|}} 
 \overline{B}( r/q,|q|^{-(1+\tau)} ). 
\end{align*}
Then the $s$-dimensional Hausdorff measure of 
$E(K,B,\tau) \cap [-1/2,1/2]^n$ is 
\begin{align*}
&
\mathcal{H}_s\rbr{   E(K,B,\tau) \cap [-1/2,1/2]^n }
%&
\leq 
\sum_{\substack{ q \in \ZZ^n \\ |q| > N}} 
\sum_{\substack{r \in \ZZ^n \\ |r| \leq |q|}} 
%(2|q|^{-(1+\tau)})^s
	\rbr{ \text{diam}( \overline{B}( r/q,|q|^{-(1+\tau)} ) ) }^s
\\
&
\leq  
\sum_{\substack{ q \in \ZZ^n \\ |q| > N}} 
(2|q|+1)^n
(2|q|^{-(1+\tau)})^s
%\\
%& 
\leq 
3^{n+s} \sum_{\substack{ q \in \ZZ^n \\ |q| > N}} |q|^{n-(1+\tau)s}. 
\end{align*}
If $s > 2n/(1+\tau)$,  
%(i.e.,  $n-(1+\tau)s < -n$), 
then the last sum goes to zero as $N \to \infty$. 
\end{proof}

\section{Algebraic Number Theory: Ideals, Norms, and Bases}\label{Algebraic Number Theory: Basics section}

Let $K$ be a number field (i.e., a finite extension field of $\QQ$) of degree $n$.

%Unlike in the integers $\ZZ$ and the Gaussian integers $\ZZ[i]$, 
%the elements of the ring of integers $\ZZ_K$ may not have unique factorization into prime elements. 
%However, the 
%set of 
%ideals of $\ZZ_K$ have unique factorizations into prime ideals. 
Unlike in the integers $\ZZ$ and the Gaussian integers $\ZZ[i]$, 
it is not generally true that every element of the ring of integers $\ZZ_K$ 
can be uniquely factored into a product of prime elements. 
However, unique factorization is recovered if we consider ideals instead of elements. 
Given two ideals $\mathfrak{a}$ and $\mathfrak{b}$ of $\ZZ_K$, 
the product $\mathfrak{a}\mathfrak{b}$ is the ideal generated 
by the set $\cbr{ab : a \in \mathfrak{a}, b \in \mathfrak{b}}$. 
An ideal $\mathfrak{p}$ with $\cbr{0} \subsetneq \mathfrak{p} \subsetneq \ZZ_K$ is called prime if, 
for every $a$ and $b$ in $\ZZ_K$ such that $ab$ is in $\mathfrak{p}$, 
at least one of $a$ and $b$ is in $\mathfrak{p}$. 
Every ideal in $\ZZ_K$ can be written uniquely 
(up to the order of the factors) as a product of prime ideals.

If $\mathfrak{a}$ is an ideal of $\ZZ_K$, 
the 
%ideal 
norm of $\mathfrak{a}$ is 
$$
N(\mathfrak{a}) = |\ZZ_K / \mathfrak{a}|. 
$$
The 
%ideal 
norm is completely multiplicative: 
If $\mathfrak{a}$ and $\mathfrak{b}$ are ideals in $\ZZ_K$, 
then $N(\mathfrak{a} \mathfrak{b}) = N(\mathfrak{a}) N(\mathfrak{b})$. 
If $\mathfrak{p}$ is a prime ideal of $\ZZ_K$, 
then $N(\mathfrak{p}) = p^f$ 
where $p$ is the unique rational prime contained in $\mathfrak{p}$  
and $f$ is the positive integer equal to the degree of 
$\ZZ_K / \mathfrak{p}$ over $\ZZ / p\ZZ$.

%For principal ideals, we have an alternate expression for the norm. 
%It requires some preliminaries. 
For principal ideals, we have an alternate expression for the norm in terms of the embeddings 
(i.e., injective homomorphisms) 
of $K$ into $\mathbb{C}$. 
Let $\mathcal{E}_K$ denote the set of all such embeddings.  
%Let $\mathcal{E}_K$ denote the set of all embeddings 
%(i.e., injective homomorphisms) 
%from $K$ into $\mathbb{C}$. 
There are precisely $n = r_1 + 2 r_2$ 
of them, 
%embeddings,  
%of $K$ into $\mathbb{C}$, 
where $r_1$ is the number of real embeddings 
and $r_2$ is the number of conjugate pairs of complex embeddings. 
We denote the real embeddings by $\rho_1, \ldots, \rho_{r_1}$ 
and denote the complex embeddings by  
$\sigma_1, \overline{\sigma_1}, \ldots, \sigma_{r_2}, \overline{\sigma_{r_2}}$. 
%If $q \in \ZZ_K$ and if $\abr{q} = q\ZZ_K$ 
%is the principal ideal of $\ZZ_K$ generated by $q$, 
If $\abr{q} = q\ZZ_K$ is the principal ideal of 
$\ZZ_K$ generated by $q \in \ZZ_K$, 
then 
$$ 
N(\abr{q}) 
=
\prod_{\tau \in \mathcal{E}_K} |\tau(q)| 
=
\left( \prod_{i=1}^{r_1} |\rho_i(q)| \right) 
\left( \prod_{i=1}^{r_2} |\sigma_i(q)|^2 \right). 
$$

Given a basis $B=\cbr{\omega_1,\ldots,\omega_{n}}$ for $K$ over $\QQ$, 
we identify $\mathbb{Q}^n$ with $K$ by identifying 
$q=(q_1,\ldots,q_n) \in \mathbb{Q}^n$ with $q=\sum_{i=1}^{n} q_i \omega_i \in K$. 
%we identify $q=(q_1,\ldots,q_n) \in \mathbb{Q}^n$ with $q=\sum_{i=1}^{n} q_i b_i \in K$. 
%This identifies $\mathbb{Q}^n$ with $K$. 
Note that we use $q$ to denote both the element of $\QQ^n$ 
and the element of $K$. 
So, for example, if $q=\sum_{i=1}^{n} q_i \omega_i \in K$, then $|q| = \max_{1 \leq i \leq n} |q_i|$. 
If this identification is also a bijection between $\mathbb{Z}^n$ and $\ZZ_K$, then $B$ is called an integral basis. 
Note that if $B$ is an integral basis, then $B \subseteq \ZZ_K$.

The following simple estimates 
will be used in several places. 

\begin{lemma}\label{C1 bounds lemma}
Let $B = \cbr{\omega_1,\ldots,\omega_n}$ be a basis for $K$ over $\QQ$. Define 
\begin{align*}%\label{defn of C1}
C_B = \max\cbr{ \textstyle \sum_{i=1}^{n} |\tau(\omega_i)| : \tau \in \mathcal{E}_K }.
\end{align*}
Let $q \in \ZZ_K$. 
For every $\tau \in \mathcal{E}_K$,
\begin{align}\label{C1 bound}
|\tau(q)| 
\leq C_B |q|.   
\end{align}
Consequently, 
\begin{align}\label{C1 norm bound}
N(\abr{q}) \leq C_B^n |q|^n.
\end{align}
\end{lemma}
\begin{proof}
If $q = \sum_{i=1}^{n} q_i \omega_i \in K$, then 
$\tau(q) = \sum_{i=1}^{n} q_i \tau(\omega_i) \in \CC$. 
Thus $|\tau(q)| \leq \max_{1 \leq i \leq n} |q_i| \sum_{i=1}^n |\tau(\omega_i)| = C_B|q|$. 
\end{proof}

\section{Algebraic Number Theory: A Divisor Bound}
\label{Algebraic Number Theory: A Divisor Bound section}

The main result of this section is 
%the following divisor bound, 
the divisor bound Proposition \ref{divisor bound thm}, 
which we need for the proof of Theorem \ref{main-thm-2}. 
It is essentially due to Elkies \cite{ElkiesMO11}. 

Let $K$ be a number field of degree $n$ over $\QQ$. 
%Define $r = r_1 + r_2 - 1$, 
%where $r_1$ and $r_2$ are (respectively) 
%the number of real embeddings and 
%the number of pairs of complex conjugate embeddings 
%of $K$ into $\CC$. 
For $s,t > 0$, define 
\begin{align*}
w_s(t) = \exp\rbr{  \frac{s \log t}{\log \log t} }. 
\end{align*}

\begin{prop}\label{divisor bound thm} 
Let $B=\cbr{\omega_1,\ldots,\omega_{n}}$ be a basis for $K$.  
For every $s \in \ZZ_K$ and $M \geq 2$, define 
%$D(s,M)$ to be the set of $q \in \ZZ_K$ 
%such that $q$ divides $s$ and $|q| \leq M$. 
$$
D(M,s) = \cbr{q \in \ZZ_K : q \mid s \text{ and } |q| \leq M}. 
$$
For every $s \in \ZZ_K$, $M \geq 2$, and $\zeta > \log 2$, 
$$
|D(M,s)| 
%\leq C_K 
\lesssim_{B,\zeta} 
w_{\zeta}( N(s) ) \log^{r_1+r_2-1} (M).  
$$
\end{prop}

The rest of this section is devoted to the proof of Proposition \ref{divisor bound thm}. 
The key idea is that $q$ divides $s$ in $\ZZ_K$ 
if and only if $q$ generates a principal ideal that divides the ideal $\abr{s}$. 
%So, to bound $|D(s,M)|$, we just need to bound two things: 
Therefore, we bound $|D(M,s)|$ by the product of upper bounds for the following two quantities: 
\begin{enumerate}[(i)]
\item The number of $q \in \ZZ_K$ with $|q| \leq M$ that generate a given principal ideal $\abr{a}$.  
\item The number of principal ideals $\abr{a}$ that divide $\abr{s}$. 
\end{enumerate}
%Then $|D(s,M)|$ will be bounded by the product of the bounds for (1) and (2). 
%In fact, instead of bounding (ii), we bound the larger quantity: 
%\begin{description}[font=\normalfont]
%\item[(ii)] The number of ideals that divide $\abr{s}$.
%\end{description} 

The following lemma gives us 
%an appropriate bound on (i). 
the desired upper bound on (i). 

\begin{lemma}\label{generator bound} 
Let $B=\cbr{\omega_1,\ldots,\omega_{n}}$ be a basis for $K$. 
Let $M \geq 2$. 
Let $a \in \ZZ_K$. 
Define $G(a,M)$ to be the set of elements $q \in \ZZ_K$ 
such that $|q| \leq M$ and $q$ generates $\abr{a}$ (i.e., $\abr{q} = \abr{a}$).   
Then 
%the number of elements in $G(a,M)$ is 
%$G(a,M)$ is a finite set of size 
\begin{align*}
|G(a,M)| 
%%%\leq C_K
\lesssim_B    
\log^{r_1+r_2-1}(M).
\end{align*}
%where $C_{K}$ is a constant depending only on $K$. 
\end{lemma}
\begin{proof}
%The proof is similar to the proof of Dirichlet's unit theorem. 
Define the map $\lambda: \ZZ_K \to \mathbb{R}^{r_1 + r_2}$ by 
$$
\lambda(x) = 
\left( \log|\rho_1(x)|, \ldots, |\rho_{r_1}(x)|, 2\log|\sigma_1(x)|, 
\ldots, 2\log|\sigma_{r_2}(x)| \right). 
$$
%Note $\log |N(x)|$ is the sum of the components of $\lambda(x)$. 
In the standard proof of Dirichlet's unit theorem 
(see, for example, \cite[Ch.7]{Jarvis14}), 
it is established  
that $\lambda$ is a group homomorphism,  
that the kernel of $\lambda$ 
is the finite cyclic group of roots of unity in $K$, 
and 
that $\lambda$ sends the unit group $\ZZ_K^{\times}$ 
to a $(r_1+r_2-1)$-dimensional lattice $L \subseteq \mathbb{R}^{r_1 + r_2}$.  
Since every generator of $\langle a \rangle$ in $\ZZ_K$
is the product of $a$ with some unit, 
it follows that $\lambda$ maps the 
set of generators of $\langle a \rangle$ to a translate of $L$, 
namely $\lambda(a)+L$. 
Now consider an arbitrary 
$q \in G({a,M})$ 
and $\tau \in \mathcal{E}_K$. 
By \eqref{C1 bound}, 
\begin{align*}
1 \leq N(\abr{a}) = N(\abr{q}) 
= |\tau(q)| 
\prod_{ \tau' \in \mathcal{E}_K \setminus \cbr{\tau}} | \tau'(q) | 
\leq 
|\tau(q)| (C_B M)^{n-1}. 
\end{align*}
Rearranging and taking the logarithm gives 
$-\log |\tau(q)| \leq (n-1) \log(C_B M)$. 
Also by \eqref{C1 bound}, 
$\log |\tau(q)| \leq \log(C_B M)$.  
Therefore $|\log |\tau(q)| | \leq n \log(C_B M)$. 
It follows that $\lambda$ maps $G(a,M)$ 
into the intersection of the lattice translate $\lambda(a)+L$ 
with the cube in $\RR^{r_1+r_2}$ 
centered at the origin 
and having side length $\leq 4 n \log(C_B M)$. 
Since the $(r_1+r_2-1)$-dimensional lattice $L$ 
does not depend on $\langle a \rangle$, 
the number of points of $\lambda(a)+L$ in this cube 
is $\leq C (4 n \log(C_B M))^{r_1+r_2-1}$, 
where $C$ is a constant depending only on $K$ 
(via the geometry of the lattice $L$).  
Since the kernel of $\lambda$ is finite, 
%Therefore 
the number of elements of $G(a,M)$ satisfies the same bound, 
with a larger constant $C$.  
Since $M \geq 2$, by further increasing $C$ 
relative to $r_1$,$r_2$, $n$, and $C_B$, 
we obtain a bound of the desired shape. 
\end{proof}

Now we turn to finding an appropriate upper bound on (ii). 
We introduce 
notation for 
two divisor functions. 
For every positive integer $\ell$, $d_1(\ell)$ is the number of positive integers that divide $\ell$. 
For every ideal $\mathfrak{s}$ of $\ZZ_K$, 
$d_2(\mathfrak{s})$ is the number of ideals that divide $\mathfrak{s}$ in $\ZZ_K$.

The function $d_1$ obeys the following classic divisor bound, 
which is due to Wigert. 
%For the proof, we refer the reader to 
%Wigert \cite{Wigert} or Hardy and Wright \cite[Theorem 317]{HW}. 
For a proof and historical discussion, 
we refer the reader to  Hardy and Wright \cite[Theorem 317]{HW}. 
\begin{lemma}\label{integer divisor bound}
%For every $\zeta > \log 2$, 
%there exists $L_{\zeta} \in \NN$ 
%such that, for every $\ell \in \NN$,   
For every $\ell \in \NN$ and $\zeta > \log 2$, 
$$
d_1(\ell) 
%\leq \exp\rbr{ \frac{\zeta \log \ell}{ \log \log \ell }}. 
\lesssim_{\zeta} w_{\zeta}(\ell). 
$$
\end{lemma}

The next lemma relates $d_1$ to $d_2$.

\begin{lemma}\label{ideal divisor bound}
For every ideal $\mathfrak{s}$ in $\ZZ_K$, 
$$
d_2(\mathfrak{s}) \leq d_1(N(\mathfrak{s})). 
$$
%that is, the number of ideal divisors of $\mathfrak{s}$ in $\ZZ_K$ 
%is less than or equal to 
%the number of positive integer divisors of the integer $N(\mathfrak{s})$.  
\end{lemma}
\begin{proof}
Suppose the unique factorization of $\mathfrak{s}$ 
into prime ideals in $\ZZ_K$ is 
$$
\mathfrak{s} = \mathfrak{p}_1^{e_1} \cdots \mathfrak{p}_k^{e_k}.  
$$
For each $i$, 
$N(\mathfrak{p}_i) = p_i^{f_i}$,  
where $p_i$ is a rational prime and $f_i$ is a positive integer. 
So the unique factorization of $N(\mathfrak{s})$ into rational primes is 
$$
N(\mathfrak{s}) 
= N(\mathfrak{p}_1)^{e_1} \cdots N(\mathfrak{p}_k)^{e_k} 
= p_1^{e_1 f_1} \cdots  p_k^{e_k f_k}. 
$$
Therefore 
$$
d_2(\mathfrak{s}) = \prod_{i=1}^{k} (e_i + 1) \leq \prod_{i=1}^{k} (e_i f_i + 1) = d_1( N(\mathfrak{s}) ). 
$$
\end{proof}

By combining Lemma \ref{integer divisor bound} 
with Lemma \ref{ideal divisor bound} 
and by taking $\mathfrak{s} = \abr{s}$, 
we obtain the desired upper bound on (ii). 
(In fact, we get an upper bound on a quantity larger than (ii); namely, the number of ideals (principal and otherwise) that divide $\abr{s}$.) 
This completes the proof of Proposition \ref{divisor bound thm}.

\section{Algebraic Number Theory : An Exponential Sum}
\label{Algebraic Number Theory : An Exponential Sum section}

In this section we prove an exponential sum estimate, 
namely Proposition \ref{GeomSeriesProp}, 
that we will need in the proof of Theorem \ref{main-thm-2}. 

Let $K$ be a number field of degree $n$ over $\QQ$. 
An element $q$ of $K$ induces a linear map on $K$ given by $x \mapsto qx$. 
Given a basis $B$ for $K$ over $\QQ$, we write $A_{q, B}$ 
for the matrix representation of this linear map in the basis $B$. 
The transpose of a matrix $A$ is denoted by $A^T$. 
Vectors $x,y \in \RR^n$ are viewed as 
$n \times 1$ column matrices, 
so 
%the dot product is 
$x \cdot y = x^T y$. 

The following lemma provides a formula for the exponential sum we are interested in. 
\begin{lemma}\label{GeomSeriesLemma}
Let $B = \cbr{\omega_1,\ldots,\omega_n}$ be an integral basis for $K$ over $\QQ$. 
Let $s,q \in \ZZ^n$. 
Let $R_q \subseteq \ZZ_K$ be any complete set of representatives 
of $\ZZ_K / \abr{q}$. 
Then 
\begin{align*}
\sum_{r \in R_q} e(s \cdot (r/q)) 
=
\sum_{r \in R_q} e\rbr{ ( (A_{q^{-1},B})^{T} s )^T  r } 
=
\left\{
\begin{array}{ll}
N(\abr{q}) & \text{if }  (A_{q^{-1},B})^{T} s \in \ZZ^n \\
 0       & \text{if }  (A_{q^{-1},B})^{T} s \notin \ZZ^n
\end{array}
\right.
\end{align*}
\end{lemma}
\begin{proof}
For every $r \in \ZZ^n$, 
$$
s \cdot (r/q) 
= s \cdot (A_{q^{-1},B}) r 
= s^T (A_{q^{-1},B}) r 
= ( (A_{q^{-1},B})^T s )^T r. 
$$
%Therefore 
%\begin{align}\label{GeomSeriesLemma 1}
%\sum_{r \in R_q} e(s \cdot r/q) 
%=
%\sum_{r \in R_q} e\rbr{ ( (A_{q,B}^{-1})^{T}s )^T  r }.  
%\end{align}
If $(A_{q^{-1},B})^T s \in \ZZ^n$, 
then 
$$
\sum_{r \in R_q} e(s \cdot (r/q)) 
= \sum_{r \in R_q} 1 =  |R_q| = |\ZZ_K / \abr{q}| = N(\abr{q}).  %%% = |N(q)|. 
$$
Now assume $(A_{q^{-1},B})^T s \notin \ZZ^n$. 
So some component, say the $j$-th component, of 
$(A_{q^{-1},B})^T s$ is not an integer. 
Then  $s \cdot (\omega_j/q) = (A_{q^{-1},B})^T s)^T \omega_j \notin \ZZ$,  
and so 
$e\rbr{s \cdot (\omega_j/q)} = e\rbr{ ( (A_{q^{-1},B})^T s )^T \omega_j } \neq 1$. 
%Since $B \subseteq \ZZ_K$, 
%$\omega_j \in \ZZ_K$. 
%Note that 
Now observe that 
%Therefore 
$R_q+\omega_j$  
is also a complete set of representatives of 
$\ZZ_K / \abr{q}$. %(this uses that $\omega_j$ is an algebraic integer)
So there is a bijection $\rho:R_q \to R_q + \omega_j$ 
such that, for each $r \in R_q$, there exists $k_r \in \ZZ_K$ 
such that $r=\rho(r) + k_r q$. 
Therefore 
\begin{align*}
\sum_{r \in R_q} e\rbr{ s \cdot (r/q) }
&=
\sum_{r \in R_q} e\rbr{ s \cdot ((\rho(r) + k_r q)/q) } 
=
\sum_{r \in R_q} e\rbr{ s \cdot (\rho(r)/q) }
\\
&=
\sum_{r \in R_q + \omega_j} e\rbr{ s \cdot (r/q) }
=
\sum_{r \in R_q} e\rbr{ s \cdot ((r+\omega_j)/q) }
\\
&=
e\rbr{ s \cdot (\omega_j/q) } \sum_{r \in R_q} e\rbr{ s \cdot (r/q) }
\end{align*}
Since $e\rbr{s \cdot (\omega_j/q)}  \neq 1$, 
the sum 
%appearing first and last 
must equal zero. 
\end{proof}

We know that $A_{q^{-1},B}$ is the matrix representation 
of $q^{-1}$ with respect to the basis $B$. 
To use Lemma \ref{GeomSeriesLemma}, 
we need to understand what the transpose $(A_{q^{-1},B})^T$ represents. 
We start with some examples.

\begin{example}
Let $K = \mathbb{Q}(i)$, $B = \{1, i\}$, and $q=a_0 + a_1 i$. 
Then the matrix representation of $q$ with respect to $B$ is 
\[ A_{q,B} = \left( \begin{array}{cc} 
a_0  & -a_1 \\
a_1 & a_0
\end{array} \right)\]
Notice that the transpose $(A_{q,B})^T$  is the matrix representation  
of the complex conjugate 
$\overline{q} = a_0 - a_1 i \in K$
with respect to the basis $B$. 
Note also that $(A_{q,B})^T$ is the matrix representation of 
$q$ with respect to the basis $B' = \cbr{1,-i}$. 
In other words, $(A_{q,B})^T = A_{q,B'}$. 
\end{example}

\begin{example}
Let $\omega = \sqrt[4]{-1} 
%= \frac{\sqrt{2} + \sqrt{2}i}{2}$. 
= ({\sqrt{2} + \sqrt{2}i})/{2}$. 
Let $K = \mathbb{Q}(\omega)$, $B = \{1, \omega, \omega^2, \omega^3\}$, 
and $q=a_0 + a_1\omega + a_2 \omega^2 + a_3 \omega^3$. 
Then the matrix represenation of  
$q$ with respect to $B$ is 
 \[ A_{q,B} = \left( \begin{array}{cccc} 
a_0 & -a_3 & -a_2 & -a_1 \\
a_1 &  a_0 & -a_3 & -a_2 \\
a_2 &  a_1 &  a_0 & -a_3 \\
a_3 &  a_2 &  a_1 &  a_0
\end{array} \right)\]
The transpose $(A_{q,B})^T$ is the matrix 
representation of the Galois conjugate 
$q' = a_0 - a_3 \omega - a_2 \omega^2 - a_1 \omega^3 \in K$ 
with respect to the basis $B$. 
In fact, $q'$ is the complex conjugate of $q$. 
Note also that 
$(A_{q,B})^T$ is the matrix representation of $q$ 
with respect to the basis $B' = \cbr{1,-\omega^3,-\omega^2,-\omega}$. 
In other words, $(A_{q,B})^T = A_{q,B'}$. 
\end{example}

\begin{example}
Let $K = \mathbb{Q}(\sqrt[3]{2})$, $B = \{1, \sqrt[3]{2}, \sqrt[3]{4}\}$, 
and $q = a_0 + a_1 \sqrt[3]{2} + a_2 \sqrt[3]{4}$. 
The matrix representation of $q$ with respect to $B$ is 
\[A_{q,B} = \left( \begin{array}{ccc} 
a_0 & 2 a_2 & 2 a_1\\
a_1 & a_0   & 2 a_2 \\
a_2 & a_1   & a_0 \end{array} \right)\]
%If $a_1 \neq 0$ or $a_2 \neq 0$, 
%then the transpose $(A_{q,B})^T$ is not the matrix representation of any element of 
%$K$  
%with respect to the basis $B$.  
If $a_1 \neq 0$ or $a_2 \neq 0$, then 
$(A_{q,B})^T$ is not the matrix representation 
of a conjugate of $q$ with respect to the basis $B$, 
in contrast to the previous examples. 
In fact, $(A_{q,B})^T$ is not the matrix representation of any element of 
$K$ with respect to the basis $B$. 
However, as in the previous examples, $(A_{q,B})^T$ is the matrix representation of $q$ 
with respect to the basis $B' = \cbr{1,{1}/{\sqrt[3]{2}},{1}/{\sqrt[3]{4}}}$. 
Notice that $B'$ is not an integral basis even though $B$ is. 
\end{example}

\begin{comment}
%The last example illustrates 
The examples above illustrate 
the second 
%obstacle 
point 
discussed at the end of 
Section \ref{background and motivation}. 
However, it is not completely disastrous, and, 
as the following lemma shows
things cannot get any worse. 
\end{comment}

%The next lemma shows the conclusion of the examples 
%above is actually the rule. 

The following lemma verifies the property of transposes of matrix representations that is suggested by the examples. 

\begin{lemma}\label{DualRepresentation}
%Let $B$ be a basis for $K$. 
Given any basis $B$ for $K$, 
there exists a basis $B'$ for $K$ 
such that, for every $q \in K$, 
%if $A_{q,B}$ is the matrix representation for $q$ with respect to $B$, 
%then 
$(A_{q,B})^T$ is the matrix representation for $q$ with respect to $B'$, 
i.e.,  
$(A_{q,B})^T = A_{q,B'}$. 
\end{lemma}
A statement of this lemma can be found, for example, in \cite{Faddeev78}. 
It can be obtained as a corollary of the Skolem-Noether theorem, 
but we include a proof for completeness. 

\begin{proof}[Proof of Lemma \ref{DualRepresentation}]
It is easy to check that the map $\phi_1$: $q \mapsto A_q$ is an injective ring homomorphism
from $K$ into $M_n(\QQ)$ (the ring of $n \times n$ matrices with rational entries). 
The same is true of the map $\phi_2$: $q \mapsto A_q^T$. 
It suffices to show that there exists a matrix $A$ with rational entries such that $\phi_2(x) = A^{-1} \phi_1(x) A$ for all $x \in K$. 
Note $K$ is a simple extension of $\mathbb{Q}$, say $K = \mathbb{Q}(\theta)$. 
So the homomorphisms $\phi_1$ and $\phi_2$ are entirely determined by $\phi_1(\theta)$ and $\phi_2(\theta)$, respectively. 
So it is enough to show that $\phi_1(\theta)$ and $\phi_2(\theta)$ are similar matrices. 
%By the Cayley-Hamilton theorem, $\phi_1(\theta)$ (resp. $\phi_2(\theta)$) 
%satisfies its own characteristic polynomial. 
Let $P_{1}$ (resp. $P_{2}$) 
be the characteristic polynomial of $\phi_1(\theta)$ 
(resp. $\phi_2(\theta)$). 
By the Cayley-Hamilton theorem and the fact that $\phi_1$ is a homorphism, 
%Then 
$0 = P_{1}(\phi_1(\theta)) = \phi_1(P_{1}(\theta))$.  
But, since $\phi_1$ is injective, 
it follows that $P_{1}(\theta) = 0$. 
The polynomial $P_{1}$ is a monic polynomial of degree $n$, 
and $\theta$ is an element of $K$ of degree $n$, 
so it follows that $P_{1}$ is the minimal polynomial of $\theta$. 
The same, of course, can be said for $P_{2}$, and 
thus $P_{1} = P_{2}$. 
Furthermore, $K/\mathbb{Q}$ is separable. 
So $P_{1}$ does not have any multiple roots. 
This implies $\phi_1(\theta)$ has $n$ distinct eigenvalues 
and is therefore diagonalizable. 
Thus $\phi_1(\theta)$ and $\phi_2(\theta)$ 
are diagonalizable matrices with the same eigenvalues and are therefore similar. 
A standard argument (see e.g. \cite{HornJohnson90}, Section 3.4, Exercise 3) involving the rational canonical form shows 
that the similarity matrix $A$ can be taken to be rational.
\end{proof}

For the definition of the constant $\beta$ in Proposition \ref{GeomSeriesProp} below, 
we note the following standard fact. 

\begin{lemma}
For every algebraic number $\alpha$, there is a rational integer $d$ such that $d\alpha$ is an algebraic integer. 
\end{lemma}
\begin{proof}
If $\alpha$ is a root of the polynomial 
$x^m + (a_{m-1}/b_{m-1})x^{m-1} + \cdots  + (a_0/b_0)$, where $a_i,b_i \in \ZZ$, 
then, with $d=b_{m-1} \cdots b_0$, it follows that $d\alpha$ is a root of the polynomial 
$x^m + d (a_{m-1}/b_{m-1})x^{m-1} + \cdots  + d^{m} (a_0/b_0)$, 
whose coefficients are integers. 
\end{proof}

%For the definition of the constant $D$ in Proposition \ref{GeomSeriesProp} below, 
%we note the following standard fact. 

%\begin{lemma}
%For every algebraic number $\alpha$, 
%there is a rational integer $d$ such that $d\alpha$ 
%is an algebraic integer. 
%\end{lemma}
%\begin{proof}
%If $\alpha$ is a root of the polynomial 
%$x^m + (c_{m-1}/c_{m-1})x^{m-1} + \cdots  + (c_0/d_0)$, 
%where $c_i,d_i \in \ZZ$, 
%then with $d=d_{m-1} \cdots d_0$ 
%we have that $d\alpha$ is a root of the polynomial 
%$x^m + k (c_{m-1}/d_{m-1})x^{m-1} + \cdots  + k^{m} (c_0/d_0)$, 
%whose coefficients are integers. 
%\end{proof}

Finally, we are ready for the main result of this section. 

\begin{prop}\label{GeomSeriesProp}
Let $B = \cbr{\omega_1,\ldots,\omega_n}$ be an integral basis for $K$ over $\QQ$.  
%such that $B \subseteq \ZZ_K$. 
Let $B' =  \{\omega_1', \omega_2', \ldots, \omega_{n}'\}$ 
be the basis corresponding to $B$ given by Lemma \ref{DualRepresentation}.  
Let $\beta$ be the smallest positive integer such that 
$\beta \omega_i' \in \ZZ_K$ for all $1 \leq i \leq n$. 
%Identify $K$ with $\QQ^n$ via $B$. 
Let $s \in \ZZ^n$ and define $s' = \sum_{i=1}^{n} s_i \omega_i'$. 
Let $q \in \ZZ^n$ and let $R_q \subseteq \ZZ_K$ 
be any complete set of representatives of $\ZZ_K / \abr{q}$. 
Then 
\begin{align*}
\abs{\sum_{r \in R_q} e(s \cdot r/q)} 
\leq 
\left\{
\begin{array}{ll}
N(\abr{q}) & \text{if $q \mid \beta s'$}  \\
 0       & \text{if $q \nmid \beta s'$} 
\end{array}
\right.
\end{align*}
\end{prop}
\begin{proof}
We always have 
$$
\abs{\sum_{r \in R_q} e(s \cdot r/q)} \leq |R_q| = |\ZZ_K / \abr{q}| = N(\abr{q}). %%% = |N(q)|. 
$$
We will show that $(A_{q^{-1},B})^{T}s \in \ZZ^n$ implies $q \mid \beta s'$. 
Combining the contrapositive with Lemma \ref{GeomSeriesLemma} 
will then complete the proof. 
Assume $(A_{q^{-1},B})^{T} s \in \ZZ^n$. 
By Lemma \ref{DualRepresentation}, 
this assumption is equivalent to the statement that $(A_{q^{-1},B'}) s \in \ZZ^n$. 
%Using the identification of $K$ with $\QQ^n$ via the basis $B'$, 
The last statement is equivalent to saying that 
$s'/q = a_1 \omega_1' + \cdots + a_n \omega_n'$ for some $a_1,\ldots,a_n \in \ZZ$. 
Multiplying by $\beta$ gives 
$\beta s'/q = a_1 (\beta \omega_1') + \cdots + a_n (\beta \omega_n')$. 
Since $\beta \omega_i' \in \ZZ_K$ for all $1 \leq i \leq n$, 
we have $\beta s'/q \in \ZZ_K$, i.e., $q \mid \beta s'$ in $\ZZ_K$. 
\end{proof}

\section{
\texorpdfstring{
Proof of Theorem \ref{main-thm-2}: The Sets $Q(M)$, $Q'(M)$, and $Q''(M)$
}{
Proof of Theorem \ref{main-thm-2}: The Sets Q(M), Q'(M), and Q''(M)
}
}
\label{The Sets $Q(M)$, $Q'(M)$, and $Q''(M)$ section}

We now begin the proof of Theorem \ref{main-thm-2} proper. 
Fix $\tau > 1$. 
Fix a number field $K$ of degree $n$ over $\QQ$. 
Fix an integral basis $B = \cbr{\omega_1,\ldots,\omega_n}$ 
for $K$.

From Proposition \ref{GeomSeriesProp},  
recall the definition of the basis $B'$, the constant $\beta \in \NN$,  
and the element $s' \in K$ for $s \in \ZZ^n$. 
Note that $\beta s' \in \ZZ_K$ for every $s \in \ZZ^n$. 

Let $M \geq 1$ be an arbitrary real number. 
Define 
$$
Q(M) = \cbr{q \in \ZZ^n :  M/2 < |q| \leq M}. 
$$
%Note $Q(M)$ contains $\approx M^n$ elements. 
Note 
\begin{align}\label{Q(M) card}
|Q(M)| \geq 2^{n-2} M^n. 
\end{align}

We define a new set $Q'(M)$ by removing from $Q(M)$ those $q$ which divide 
$\beta s'$ for some small non-zero $s \in \ZZ^n$.  
This is needed for Lemma \ref{FM hat lemma 3} below. 
We also show that this requires removing only a small number of elements, 
which is important for Lemma \ref{FM hat lemma 4} below. 
Define 
$$
S(M) = \cbr{s \in \ZZ^n: 0< |s| \leq M^{1/(2n)}}. 
$$
With $D(M, \beta s')$ defined as in 
Proposition \ref{divisor bound thm}, 
define 
$$
Q'(M) = Q(M) \setminus \bigcup_{s \in S(M)} D(M, \beta s').  
$$ 
Note 
\begin{align}\label{S(M) card}
|S(M)| \leq 2^n M^{1/2}. 
\end{align}
For each $s \in S(M)$, 
\eqref{C1 norm bound} gives 
$$
|N(\abr{\beta s'})| 
\leq C_{B'}^n \beta^n |s|^n 
\leq C_{B'}^n \beta^n M^{1/2},  
$$
and so Proposition \ref{divisor bound thm} implies 
\begin{align}\label{D(M,beta s') card}
|D(M,\beta s')|
= M^{o(1)}.  
\end{align}
%
%
%
%%%Then $Q'(M)$ contains $\approx M^n - M^{1/2+ o(1)} \approx M^n$ 
%%%elements for all sufficiently large $M$. 
By \eqref{Q(M) card}, \eqref{S(M) card}, and \eqref{D(M,beta s') card}, 
there is a number $M_0'$ such that, for all $M \geq M_0'$,  
\begin{align}\label{Q'(M) card 2} 
|Q'(M)| \geq 2^{n-2} M^n - 2^n M^{1/2 + o(1)} \gtrsim M^n. 
\end{align}
%for all $M \geq M_0'$. 

Now we choose a subset $Q''(M)$ of $Q'(M)$ 
consisting of elements $q$ which all 
have approximately the same norm $N(\abr{q})$. 
We also ensure that $Q''(M)$ is not too much smaller than $Q'(M)$. 
This is needed for Lemma \ref{FM hat lemma 4} below. 

By \eqref{C1 norm bound}, for every $q \in Q(M)$ 
and hence for every $q \in Q'(M)$,
 we have 
$$
1 \leq N( \abr{q}  ) 
\leq C_B^n |q|^n 
\leq C_B^n M^{n}. 
$$
Define 
\begin{align}\label{J defn}
J = \lceil \log_2( C_B^n M^{n} ) \rceil. 
\end{align}
Partition $Q'(M)$ dyadically  as 
$$
Q'(M) = \bigcup_{j=0}^{J} 
\cbr{q \in Q'(M) : 2^{-j-1}  < N(\abr{q}) C_B^{-n} M^{-n}  \leq 2^{-j}    }
$$
By the pigeonhole principle, there exists a 
$j_0(M) \in \cbr{0, \ldots, J}$ 
%such that, for all sufficiently large $M$, the set
such that the set 
$$
Q''(M) = \cbr{q \in Q'(M) : 2^{-j_0(M)-1}  < N(\abr{q}) C_B^{-n} M^{-n} \leq 2^{-j_0(M)}  }
$$ 
%contains $\gtrsim |Q(M)| / K \approx M^n / \log M$ elements. 
has cardinality 
$
|Q''(M)| \geq  |Q'(M)| / J. 
$
Therefore, by \eqref{Q'(M) card 2} and \eqref{J defn}, 
there is a number $M_0''$ such that, for all $M \geq M_0''$,   
\begin{align}\label{Q''(M) card 2}
|Q''(M)| \gtrsim \frac{M^n}{\log M}. 
\end{align}
%for all $M \geq M_0''$. 

\section{
\texorpdfstring{
Proof of Theorem \ref{main-thm-2}: The Function $F_M$
}
{
Proof of Theorem \ref{main-thm-2}: The Function FM
}
}
\label{The Function FM section}

Fix $\phi: \RR^n \to \RR$ such that $\phi$ is $C^{\infty}$, $\phi \geq 0$, $\int \phi = 1$, $\supp(\phi) \subseteq [-1,1]^n$.  
%For $\epsilon > 0$, define $\phi_{\epsilon}(x) = \epsilon^{-d} \phi(x/\epsilon)$. 
%Let $M > 0$.  
%Define 
%$$
%Q(M) = \cbr{q \in \ZZ^n :  M/2 < |q| \leq M}
%$$
%Note $Q(M)$ contains $\approx M^n$ elements. 
Define $c_M$ by 
$$
\frac{1}{c_M} =  \sum_{q \in Q''(M)} |N(q)|. 
$$
Define $\epsilon_M = M^{-(1+\tau)}$ and 
\begin{align}\label{FM defn}
F_M(x) = c_M \sum_{q \in Q''(M)} \sum_{r \in \ZZ^n} 
\epsilon_M^{-n} \phi ((x-r/q)/\epsilon_M)  
\end{align}
for each $x \in \RR^n$. 
%Here $c_M$ is a constant to be specified later. 
Notice that $\epsilon_M^{-n} \phi ((x-r/q)/\epsilon_M) $ is 
an $L^1$-normalized bump function on the 
$\ell^{\infty}$-ball with radius $\epsilon_M = M^{-(1+\tau)}$ and center $r/q$. 
Observe that $F_M$ is $\ZZ^n$-periodic, $F_M$ is $C^{\infty}$, and $F_M \geq 0$. 
Note also that, for each fixed $q$, 
the inner sum in the definition of $F_M$ has only finitely many non-zero terms 
because $\supp(\phi) \subseteq [-1,1]^n$. 

\begin{lemma}\label{support lemma}
For all $M \geq 1$, 
\begin{align}
\label{2.6}
\supp(F_M) \subseteq 
\bigcup_{q \in Q''(M)} \bigcup_{r \in \ZZ^n} 
\{ x \in \RR^n :  |x - r/q| \leq |q|^{-(1+\tau)}   \}. 
\end{align}
For any sequence $(M_k)_{k=1}^{\infty}$ with $2M_k \leq M_{k+1}$ for all $k \in \NN$, 
\begin{align}
\label{2.7}
\bigcap_{k=1}^{\infty} \supp(F_{M_k}) \subseteq E(K,B,\tau). 
\end{align}
\end{lemma}
\begin{proof}
%The set on the right of \eqref{2.6} is closed. Thus, to prove \eqref{2.6}, it suffices to prove it with $\supp(F_M)$ replaced by $\cbr{x \in \RR^n : F_M(x) > 0}$
Let $x \in \RR^n$. 
Since $\phi \geq 0$ and $\supp(\phi) \subseteq [-1,1]^n$, 
if $F_M(x) > 0$, then there exist $q \in Q''(M)$ and $r \in \ZZ^n$ such that 
$|\epsilon_M^{-1}(x-r/q)| \leq 1$, and 
hence $|x-r/q| \leq \epsilon_M = M^{-(1+\tau)} \leq |q|^{-(1+\tau)}$. 
This proves \eqref{2.6} with $\cbr{x \in \RR^n : F_M(x) > 0}$ in place of $\supp(F_M)$. 
But, since the set on the right of \eqref{2.6} is closed, this actually proves \eqref{2.6}. 
%But \eqref{2.6} follows immediately because the set on the right of \eqref{2.6} is closed. 
If $x \in \supp(F_{M_k})$ for every $k \in \NN$, 
then for every $k \in \NN$ we get a pair $(q_k,r_k) \in Q''(M_k) \times \ZZ^n$ with 
$|x - r_k/q_k| \leq |q_k|^{-(1+\tau)}$. The pairs must be distinct because 
$$
|q_k| \leq M_k \leq M_{k+1}/2 < |q_{k+1}|
$$
for all $k \in \NN$. This proves \eqref{2.7}. 
\end{proof}

\section{
\texorpdfstring{
Proof of Theorem \ref{main-thm-2}: The Fourier Transform of $F_M$
}
{
Proof of Theorem \ref{main-thm-2}: The Fourier Transform of FM
}
}
\label{The Fourier Transform of FM section}

%Recall that, in the definition of $E(K,B,\tau)$, 
%Recall that 
%$B$ is an integral basis for $K$ over $\QQ$.  
%Thus $\ZZ_K$ is identified with $\ZZ^n$. 
For each $q \in \ZZ^n$, let $R_q$ be a fixed set of representatives 
of $\ZZ_K / \abr{q}$. 
Note that the cardinality of $R_q$ is 
\begin{align}\label{Rq card}
|R_q| = | \ZZ_K / \abr{q} | = N(\abr{q}). %%%= |N(q)|. 
\end{align}

\begin{lemma}\label{FM hat lemma}

For all $M \geq 1$ and $s \in \ZZ^n$, 
\begin{align*}
\widehat{F_M}(s)
=
c_M \widehat{\phi}(s/M^{1+\tau}) \sum_{q \in Q''(M)} \sum_{r \in R_q} e(s \cdot r/q). 
\end{align*}
\end{lemma}
\begin{proof}
Since $\ZZ^n$ is identified with $\ZZ_K$ via the integral basis $B$, 
every element $r \in \ZZ^n$ can be written uniquely as $r=r'+kq$, where $r' \in R_q$ and $k \in \ZZ^n$. 
Using this and \eqref{FM defn}, we find 
\begin{align*}
\widehat{F_M}(s)
&=
c_M \sum_{q \in Q''(M)} \sum_{r \in \ZZ^n} 
\int_{[0,1]^n} \phi ((x-r/q)/\epsilon_M) e(s \cdot x) \epsilon_M^{-n}   dx
\\
&= 
c_M \sum_{q \in Q''(M)} \sum_{r \in R_q} \sum_{k \in \ZZ^n}
\int_{[0,1]^n} \phi ((x + k-r/q)/\epsilon_M) e(s \cdot x) \epsilon_M^{-n}   dx
\\
&= 
c_M \sum_{q \in Q''(M)} \sum_{r \in R_q} 
\int_{\RR^n} \phi ((x - r/q)/\epsilon_M) e(s \cdot x) \epsilon_M^{-n}   dx
\\
&= 
c_M \sum_{q \in Q''(M)} \sum_{r \in R_q} e(s \cdot r/q)
\int_{\RR^n} \phi (u)  e(\epsilon_M s \cdot u) du 
\\ 
&= 
c_M \widehat{\phi} (\epsilon_M s) \sum_{q \in Q''(M)} \sum_{r \in R_q} e(s \cdot r/q).
\end{align*}
\end{proof}

\begin{lemma}\label{FM hat lemma 1}
For all $M \geq 1$ and $s \in \ZZ^n$,  
\begin{align}
\label{F-1}
\widehat{F_M}(0) &= 1,  
%\quad \text{if } s=0 \in \ZZ^n 
\\
\label{F-2}
|\widehat{F_M}(s)| &\leq 1. 
%\quad \text{for all } s \in \ZZ^n.  
\end{align}
\end{lemma}
\begin{proof}
By \eqref{Rq card}, Lemma \ref{FM hat lemma}, and the definition of $c_M$, 
$$
\widehat{F_M}(0) 
= c_M \widehat{\phi}(0) \sum_{q \in Q''(M) } \sum_{r \in R_q} e(0 \cdot r/q) 
= c_M \sum_{q \in Q''(M) } N(\abr{q}) = 1. 
$$
Then, for every $s \in \mathbb{Z}^n$, we have 
$$
|\widehat{F_M}(s)| \leq \int_{[0,1]^n} |F_M(x) e(s \cdot x)| dx 
= \int_{[0,1]^n} F_M(x) = \widehat{F_M}(0) = 1. 
$$
%If $0 < |s| \leq M^{1/2n}$, 
%then by the definition of $Q(M)$ 
%(more specifically, by the definition of the intermediate $Q'(M)$) 
\end{proof}

\begin{comment}
\begin{lemma}\label{FM hat lemma 2}
For every $N > 0$ there exists $C_N > 0$ such that 
\begin{align*}
|\widehat{F_M}(s)|
\leq 
C_N (1+|s|)^{-N}
\end{align*}
whenever $M>0$, $s \in \ZZ^n$, $|s| \geq M^{2(1+\tau)}$. 
\end{lemma}
%
\begin{proof}
Since $\phi$ is $C^{\infty}$ and compactly supported, 
for every $N > 0$, there is a $C_N > 0$ such that 
$$
|\widehat{\phi}(\xi)| \leq C_N (1+|\xi|)^{-2N}
$$
for all $\xi \geq \RR^n$. 
%
By Lemma \ref{FM hat lemma} and the definition of $c_M$,  
\begin{align*}
|\widehat{F_M}(s)| 
&\leq  
c_M |\widehat{\phi}(s/M^{1+\tau})| 
\sum_{q \in Q(M)} \sum_{r \in R_q} |e(s \cdot r/q)| 
= 
|\widehat{\phi}(s/M^{1+\tau})| 
\\
&\leq 
C_N (1+|s|M^{-(1+\tau)})^{-2N}
\leq
C_N (1+|s|^{1/2})^{-2N} 
\leq 
C_N (1 + |s|)^{-N}
\end{align*}
\end{proof}
\end{comment}

\begin{lemma}\label{FM hat lemma AN}
For all $M \geq 1$ and $s \in \ZZ^n$, 
\begin{align*}
|\widehat{F_M}(s)|
\leq 
c_M |\widehat{\phi}(s/M^{1+\tau})| \sum_{\substack{ q \in Q''(M) \\ q \mid \beta s' } }  N(\abr{q}).  
\end{align*}
\end{lemma}
\begin{proof}
Combine Proposition \ref{GeomSeriesProp} 
and Lemma \ref{FM hat lemma}. 
\end{proof}

\begin{lemma}\label{FM hat lemma 3}
For all $M \geq 1$ and $s \in \ZZ^n$, 
if $0 < |s| \leq M^{1/2n}$,  
then 
$
\widehat{F_M}(s) = 0.   
$
\end{lemma}
\begin{proof}
%Since $Q''(M) \subsetq Q'(M)$, 
%the definition of $Q'(M)$ implies 
Because $Q''(M) \subseteq Q'(M)$ 
and because of the definition of $Q'(M)$, 
if $0 < |s| \leq M^{1/2n}$, 
then the sum over $q$ in Lemma \ref{FM hat lemma AN} 
is empty; hence, $\widehat{F_M}(s)=0$. 
\end{proof}

\begin{lemma}\label{FM hat lemma 4}
%Let $M \geq 1$ and $s \in \ZZ^n$. 
For all $\zeta > \log 2$, $M \geq \max\cbr{M_0'',2}$, and $s \in \ZZ^n$,  
\begin{align*}
|\widehat{F_M}(s)|
\lesssim_{\zeta}   
(1+|s|)^{-n/(1+\tau)} 
%|D(\beta s,M)| \log(M). 
w_{\zeta}(N(\abr{\beta s'})) 
\log^{r_1+r_2}(M). 
\end{align*}
\end{lemma}
\begin{proof}
Since $\phi$ is $C^{\infty}$ with compact support, 
$$
|\widehat{\phi}(\xi)| \lesssim (1+|\xi|)^{-n/(1+\tau)} 
$$
for all $\xi \in \RR^n$. 
Thus 
\begin{align*}%\label{lemma 4 1}
|\widehat{\phi}(s/M^{1+\tau})| 
\lesssim (1+|s|M^{-(1+\tau)})^{-n/(1+\tau)} 
\leq M^n (1+|s|)^{-n/(1+\tau)}.  
\end{align*}
%By the definition of $Q''(M)$, 
By \eqref{Q''(M) card 2} and the definition of $Q''(M)$, 
$$
\frac{1}{c_M} = \sum_{q \in Q''(M) }  N(\abr{q})
\geq 
2^{-j_0(M) - 1} C_B^n M^n  |Q''(M)| 
\gtrsim 
2^{-j_0(M) - 1} C_B^n M^n \frac{M^{n}}{\log M} 
$$
and 
$$
\sum_{\substack{ q \in Q''(M) \\ q \mid \beta s' } }  N(\abr{q}) 
\leq
2^{-j_0(M)} C_B^n M^n \sum_{\substack{ q \in Q''(M) \\ q \mid \beta s' } }  1  
\leq 
2^{-j_0(M)} C_B^n M^n |D(M,\beta s')|,  
$$
where $D(M,\beta s')$ is defined as in Proposition \ref{divisor bound thm}. 
Combining the estimates above 
with Proposition \ref{divisor bound thm} and Lemma \ref{FM hat lemma AN} gives 
%\begin{align}
%|\widehat{F_M}(s)|
%\lesssim  
%(\log M) (1+|s|)^{-n/(1+\tau)}  |D(M,\beta s)|. 
%\end{align}
%Then Proposition \ref{divisor bound thm} implies 
the desired result. 
\end{proof}

\section{Proof of Theorem \ref{main-thm-2}: Recursive Estimate}

The proposition proved in this section 
%The following proposition 
will be used recursively 
to define the measure $\mu$ in Section \ref{mu section}. 

Define 
%$$
%g(s) = (1+|s|)^{-n/(1+\tau)} \exp( \log (3+|s|)/\log \log (3+|s|))\log^{r_1+r_2+1}(3+|s|) 
%$$
$$
g(x) = 
\left\{
\begin{array}{ll}
|x|^{-n/(1+\tau)} 
%\exp( n \log |x|/\log \log |x|)
w_n(|x|)
\log^{r_1+r_2}(|x|) & \text{if } x \in \mathbb{R}^n, |x| > 3 \\
1 & \text{if } x \in \mathbb{R}^n, |x| \leq 3
\end{array} \right.
$$

\begin{prop}\label{main-lemma}
For every $\delta > 0$, $M_0 > 0$, and $\chi \in C^{\infty}_{c}(\mathbb{R}^n)$, 
there is an $M_{\ast} = M_{\ast}(\delta,M_0,\chi) \in \NN$ 
such that $M_{\ast} \geq M_0$ and 
\begin{align*}
|\widehat{\chi F_{M_{\ast}}}(\xi) - \widehat{\chi}(\xi)| \leq \delta g(\xi) \quad \text{ for all } \xi \in \mathbb{R}^{n}. 
\end{align*}
%where
\end{prop}
%
%The proof will show $M_{\ast}$ can be taken to be any 
%sufficiently large positive number.
%
\begin{proof}
We begin by recording two auxiliary estimates. 
Since $\chi \in C^{\infty}_{c}(\mathbb{R}^{n})$, for every $N > 0$, we have  
\begin{align}\label{108}
|\widehat{\chi}(\xi)| \lesssim_N (1+|\xi|)^{-N} \quad \text{ for all } \xi \in \mathbb{R}^{n}.
\end{align}
For every $p > n$, we have
\begin{align}\label{108-2}
\sup_{\xi \in \mathbb{R}^{n}} \sum_{\ell \in \ZZ^{n}} (1+|\xi - \ell|)^{-p} < \infty.
\end{align}
Fix $\xi \in \mathbb{R}^n$. 
We will write $\widehat{\chi F_{M}}(\xi) - \widehat{\chi}(\xi)$ in another form.  
Since $F_M$ is $C^{\infty}$ and $\ZZ^n$-periodic, we have
$$
F_M(x) 
= \sum_{\ell \in \ZZ^{n}} \widehat{F_M}(\ell) e({-\ell \cdot x}) 
\quad \text{ for all } x \in \mathbb{R}^{n}
$$
with uniform convergence. 
Since $\chi \in L^1(\mathbb{R}^{n})$, 
multiplying by $\chi$ and taking the Fourier transform yields 
\begin{align*}
\widehat{\chi F_M}(\xi) 
= \sum_{\ell \in \ZZ^{n}} \widehat{F_M}(\ell) \int_{\mathbb{R}^{n}} 
\chi(x) e^{-2 \pi i (\xi - \ell) \cdot x} dx
= \sum_{\ell \in \ZZ^{n}} \widehat{F_M}(\ell) \widehat{\chi}(\xi-\ell).
\end{align*} 
By Lemma \ref{FM hat lemma 1}  and \ref{FM hat lemma 3}, we have 
\begin{align}\label{110-2}
\widehat{\chi F_{M}}(\xi) - \widehat{\chi}(\xi)
=
\sum_{\ell \in \ZZ^{n}} \widehat{\chi}(\xi-\ell) \widehat{F_M}(\ell) - \widehat{\chi}(\xi) 
%%%\\
=
\sum_{|\ell| > M^{1/2n}} \widehat{\chi}(\xi-\ell) \widehat{F_M}(\ell).
\end{align}
for all sufficiently large $M$.

Fix $N > n + n/(1+\tau)$ 
and define the positive number $\eta$ by 
$N = 2\eta + n + n/(1+\tau)$. 
We estimate $\widehat{\chi F_{M}}(\xi) - \widehat{\chi}(\xi)$ 
by considering two cases.

\textbf{Case 1:} $|\xi| < \frac{1}{2} M^{1/2n}$. 

If $|\ell| > M^{1/2n}$, 
then $|\xi - \ell| \geq |\ell| - |\xi| > \frac{1}{2} M^{1/2n} > |\xi|$. 
Hence, by \eqref{F-2}, \eqref{108}, \eqref{108-2}, 
and \eqref{110-2}, we have 
\begin{align*}
&| \widehat{\chi F_{M}}(\xi) - \widehat{\chi}(\xi) |
\lesssim
\sum_{ |\ell| > M^{1/2n} } (1+|\xi - \ell|)^{-N}
= 
\sum_{ |\ell| > M^{1/2n} } 
(1+|\xi - \ell|)^{-2\eta - n - n/(1+\tau)} \\
&\leq 
(1 + |\xi|)^{-n/(1+\tau)} 
(1 + \frac{1}{2} M^{1/2n})^{-\eta} 
\sum_{|\ell| > M^{1/2n}} (1+|\xi - \ell|)^{-(n+\eta)} 
\leq
\delta g(\xi)
\end{align*}
for all sufficiently large $M$.

\textbf{Case 2:} $|\xi| \geq \frac{1}{2} M^{1/2n}$. 

Using \eqref{110-2}, write 
\begin{align*}
\widehat{\chi F_{M}}(\xi) - \widehat{\chi}(\xi) 
= 
S_1 + S_2
= 
\sum_{ \substack{  |\ell| > M^{1/2n}  \\ |\ell| \leq \frac{1}{2}|\xi|  } } \widehat{\chi}(\xi-\ell)  \widehat{F_M}(\ell) 
+
\sum_{ \substack{  |\ell| > M^{1/2n}  \\ |\ell| > \frac{1}{2}|\xi|  } } \widehat{\chi}(\xi-\ell)  \widehat{F_M}(\ell).   
\end{align*}

We first bound $S_1$. 
If $|\ell| \leq \frac{1}{2}|\xi|$, 
then $|\xi - \ell| \geq \frac{1}{2}|\xi| \geq \frac{1}{4} M^{1/2n}$. 
Hence by \eqref{F-2}, \eqref{108}, and \eqref{108-2} we have
\begin{align*}
&|S_1|
\lesssim
\sum_{\substack{ |\ell| > M^{1/2n} \\ |\ell| \leq \frac{1}{2}|\xi|} } 
(1+|\xi - \ell|)^{-N} 
= 
\sum_{\substack{ |\ell| > M^{1/2n} \\ |\ell| \leq \frac{1}{2}|\xi|}} 
(1+|\xi - \ell|)^{-2\eta - n - n/(1+\tau)} 
\\
&\leq
(1 + \frac{1}{2}|\xi|)^{-n/(1+\tau)} (1 + \frac{1}{4}  M^{1/2n} )^{-\eta} 
\sum_{ \substack{ |\ell| > M^{1/2n} \\ |\ell| \leq \frac{1}{2}|\xi|} } 
(1+|\xi - \ell|)^{-(n+\eta)} 
\leq
\frac{1}{2} \delta g(\xi)
\end{align*}
for all sufficiently large $M$.

Now we bound $S_2$. 
Fix $\zeta$ such that $\log 2 < \zeta < 1$.   
By \eqref{C1 norm bound} and Lemma \ref{FM hat lemma 4},  
%we have 
\begin{align*}
|S_2|
\lesssim 
\sum_{ \substack{  |\ell| > M^{1/2n}  \\ |\ell| > \frac{1}{2}|\xi|   }  } 
(1+|\ell| )^{-n/(1+\tau)} 
w_{\zeta}( C_{B'}^n \beta^n |\ell|^n  )
\log^{r_1+r_2}(M) 
|\hat{\chi}(\xi-\ell)|. 
\end{align*}
Note that $\log^{r_1+r_2}(x)$ is increasing and that $(1+ x )^{-n/(1+\tau)} 
w_{\zeta}( C_{B'}^n \beta^n x^n  )$ is eventually decreasing.   
Since $|\ell| > \frac{1}{2} |\xi| \geq \frac{1}{4}  M^{1/2n}$ in the sum, 
taking $M$ sufficiently large gives 
\begin{align*}
|S_2|
\lesssim 
(1+|\xi|/2 )^{-n/(1+\tau)} 
w_{\zeta}( C_{B'}^n \beta^n (|\xi|/2)^n  )
\log^{r_1+r_2}((2|\xi|)^{2n}) 
\sum_{ \substack{  |\ell| > M^{1/2n}  \\ |\ell| > \frac{1}{2}|\xi|   }  } 
|\hat{\chi}(\xi-\ell)|. 
\end{align*}
By \eqref{108} and \eqref{108-2}, 
%the sum in the last line is 
the last sum is 
$\lesssim 1$. 
Since $\zeta < 1$ and $|\xi| \geq \frac{1}{2}  M^{1/2n}$, 
taking $M$ sufficiently large gives 
$$
|S_2| \leq \frac{1}{2} \delta g(\xi). 
$$
\end{proof}

\section{
\texorpdfstring{
Proof of Theorem \ref{main-thm-2}: The Measure $\mu$
}
{
Proof of Theorem \ref{main-thm-2}: The Measure mu
}
}
\label{mu section}

In this section we construct the measure $\mu$ and prove it satisfies the desired support and Fourier decay properties.

Let $f_0:\RR^n \to \RR$ 
be a non-negative compactly supported $C^{\infty}$ function 
with $\int_{\RR^{n}} f_0(x)dx = 1$.  
With the notation of Lemma \ref{main-lemma}, define
$$
M_1 = M_{\ast}(2^{-2},1,f_0), \quad M_{k} = M_{\ast}(2^{-k-1},2M_{k-1},f_0F_{M_1}\cdots F_{M_{k-1}}) \text{ for } k = 2,3,\ldots. 
$$
Define measures $\mu_k$ on $\RR^n$ by
$$
d\mu_0 = f_0 dx, \quad d\mu_k = f_0  F_{M_1} \cdots F_{M_{k}} dx \quad \text{ for all } k \in \NN.
$$
By Lemma \ref{main-lemma}, 
\begin{align}\label{7}
|\widehat{\mu_k}(\xi) - \widehat{\mu_{k-1}}(\xi)| 
\leq 2^{-k-1} g(\xi) \quad \text{ for all } k \in \NN, \xi \in \RR^{n}. 
\end{align}
Since $g$ is bounded, \eqref{7} implies 
$(\widehat{\mu_k})_{k = 0}^{\infty}$ is a Cauchy sequence in the supremum norm.
%\eq{7} implies $(\widehat{\mu_k}(\xi))_{k \in \NN}$ is a Cauchy sequence for each fixed $\xi$.
Therefore, since each $\widehat{\mu_k}$ is a continuous function, $\displaystyle{\lim_{k \rightarrow \infty}} \widehat{\mu_k}$ is a continuous function.
By \eqref{7}, we have
\begin{align}\label{8}
|\lim_{k \rightarrow \infty} \widehat{\mu_k}(\xi) - \widehat{\mu_{\ell-1}}(\xi)| 
%&\leq
\leq 
\sum_{k={\ell}}^{\infty} |\widehat{\mu_k}(\xi) - \widehat{\mu_{k-1}}(\xi)| %\\
%\notag
%&\leq 
\leq
g(\xi) \sum_{k=\ell}^{\infty} 2^{-k-1} = 2^{-\ell} g(\xi)
\end{align}
for all $\xi \in \RR^{n}$ and $\ell \in \NN$.  
Since $\widehat{\mu_0}(0) = \int_{\RR^{n}} f_0(x)dx = 1$ and $g(0) = 1$, 
it follows from \eqref{8} that 
$$
\frac{1}{2} \leq | \displaystyle{\lim_{k \rightarrow \infty}} \widehat{\mu_k}(0) | \leq \frac{3}{2}.
$$ 
Therefore, by L\'{e}vy's continuity theorem, 
$(\mu_k)_{k=0}^{\infty}$ converges weakly 
(i.e., in distribution) to a finite non-zero Borel measure $\mu$.  
Then, by Lemma \ref{support lemma}, 
\begin{align*}
\text{supp}(\mu)  
\subseteq \bigcap_{k=0}^{\infty} \text{supp}(\mu_k)  
= \text{supp}(f_0) \cap \bigcap_{k=1}^{\infty} \text{supp}(F_{M_k}) 
\subseteq E(K,B,\tau). 
\end{align*}
Moreover, 
\begin{align*}
\widehat{\mu}(\xi) 
= \lim_{k \rightarrow \infty} \widehat{\mu_k}(\xi) 
\quad \text{ for all } \xi \in \RR^n.  
\end{align*}
Let $\epsilon > 0$ be given. 
Choose $k_{\epsilon} \in \NN$ such that $2^{-k_{\epsilon}} \leq \epsilon$. 
By \eqref{8}, we have 
\begin{align}\label{o-1}
|\widehat{\mu}(\xi) - \widehat{\mu_{k_{\epsilon}}}(\xi)|
\leq  2^{-k_{\epsilon} - 1} g(\xi) 
\leq \frac{\epsilon}{2}  g(\xi) \quad \text{ for all }  \xi \in \RR^{n}. 
\end{align}
On the other hand, since $f_0  F_{M_1} \cdots F_{M_{k_{\epsilon}}}$ is $C^{\infty}$ and compactly supported, we have 
\begin{align}\label{o-2}
|\widehat{\mu_{k_{\epsilon}}}(\xi)| \lesssim (1+|\xi|)^{-n/(1+\tau)} \quad \text{ for all } \xi \in \RR^{n}.
\end{align}
%Note that the implied constant in \eqref{o-2} depends on $\epsilon$. 
By combining \eqref{o-1} and \eqref{o-2}, we see that 
$|\widehat{\mu}(\xi)| \leq \epsilon g(\xi)$
for all sufficiently large $\xi \in \RR^{n}$, and hence 
$|\widehat{\mu}(\xi)| = o(g(\xi))$ as $|\xi| \to \infty$. 
%
%
%
%Since $f_0 \in C_{c}^{K}(\RR^{mn})$ and $K \geq a$, we have $\widehat{\mu_0}(\xi) \lesssim (1+|\xi|)^{-a}$ for all $\xi \in \RR^{mn}$. Combining this with \eq{8} gives 
%$$
%|\widehat{\mu}(\xi)| \lesssim g(\xi) \quad \forall \xi \in \RR^{mn}.
%$$
%By a standard argument (see Kahane \cite[pp.252-253]{kahane-book}), 
%it follows that $|\widehat{\mu}(\xi)| = o(g(\xi))$ as $|\xi| \to \infty$ with $\xi \in \RR^n$. 
By multiplying $\mu$ by a constant, we can make $\mu$ a probability measure. 
This completes the proof of Theorem \ref{main-thm-2}. 
%Theorems \ref{main-thm-1} and \ref{main-thm-2}.

\section{Acknowledgements}

This material is based on work supported by the National Science Foundation under Award No. 1803086.

\bibliographystyle{myplain}
\bibliography{Approximation_In_Number_Fields}
\end{document}